\documentclass[preprint]{elsarticle}

\usepackage{amsmath,amssymb,amsfonts}
\usepackage{algorithmic}
\usepackage{graphicx}
\usepackage{textcomp}
\usepackage{xcolor}
\usepackage{algorithm,algorithmic}
\usepackage{bm}
\usepackage{subfigure}
\usepackage{fullpage}
\usepackage{multirow}
\begin{document}

\begin{frontmatter}

\title{Low rank tensor completion with sparse regularization in a transformed domain}

\author[lab1]{Ping-Ping Wang}
\ead{wppunique@outlook.com}
\author[lab1]{Liang Li\corref{cor1}}
\ead{plum\_liliang@uestc.edu.cn, plum.liliang@gmail.com}
\author[lab1]{Guang-Hui Cheng}
\ead{cgh612@126.com}

\address[lab1]{School of Mathematical Sciences, University of
               Electronic Science and Technology of China,
               Chengdu, P.R. China}

\cortext[cor1]{Corresponding author}

\begin{abstract}
Tensor completion is a challenging problem with various applications. Many related models based on the low-rank prior of the tensor have been proposed. However, the low-rank prior may not be enough to recover the original tensor from the observed incomplete tensor. In this paper, we prose a tensor completion method by exploiting both the low-rank and sparse prior of tensor. Specifically, the tensor completion task can be formulated as a low-rank minimization problem with a sparse regularizer. The low-rank property is depicted by the tensor truncated nuclear norm based on tensor singular value decomposition (T-SVD) which is a better approximation of tensor tubal rank than tensor nuclear norm. While the sparse regularizer is imposed by a $\ell_{1}$-norm in a discrete cosine transformation (DCT) domain, which can better employ the local sparse property of completed data. To solve the optimization problem, we employ an alternating direction method of multipliers (ADMM) in which we only need to solve several subproblems which have closed-form solutions. Substantial experiments on real world images and videos show that the proposed method has better performances than the existing state-of-the-art methods.
\end{abstract}

\end{frontmatter}

\section{Introduction}
Estimating missing data from very limited information of observed data has attracted considerable interest recently. This problem arises from various kinds of applications in signal processing and machine learning \cite{cichocki2015tensor,zhou2015tensor,li2019efficient,Ji2016tensor}, such as image recovery, video denosing, recommender systems, and data mining. However, estimating the missing values without any prior information about the data is usually an ill-posed problem. There are many commonly adopted assumptions which can be divided into local and global information to alleviate the problem. To utilize the local information, the statistical or structural information \cite{Coupier2015image} of the observed data are used to build up the relation between the missing data and the known data, but it is obviously that the approach only focuses on local relations. It is necessary to consider the global structural information of the observed data.

In many real applications, the signals lie in a low dimensional space, for example, the natural images data have a low-rank structure \cite{hu2013fast,jiang2018matrix,zheng2019lowrank}. As a result, the matrix completion problem can be modeled as a low-rank minimization problem
\begin{equation}\label{E1}
\begin{aligned}
&\mathop{\rm{min}}\limits_{\bm{X}}{\rm{rank}}(\bm{X}) &\\
&s.t. \bm{X}_{\bm{\Omega}}=\bm{M}_{\bm{\Omega}} &\\
\end{aligned}
\end{equation}
where $\bm{X}\in\mathbb{R}^{m\times n}$, rank($\cdot$) denotes the rank of the matrix $\bm{X}$ and $\bm{\Omega}$ is the set of locations corresponding to the observed data. However, the rank function of matrix is a nonconvex and discontinuous function \cite{hu2013fast}, so the resulting (\ref{E1}) is a NP-hard problem. Theoretical studies show that the nuclear norm, i.e., the sum of singular values of a matrix, is the convex surrogate of the rank function \cite{recht2010guaranteed}. Furthermore, there are some efficient methods to solve the nuclear norm minimization problem \cite{cai2010singular}. Unfortunately, these nuclear norm methods may lead to suboptimal results, since all the singular values are treated differently when added together and minimized simultaneously while in the rank minimization process all the singular values have the same \cite{hu2013fast}. Therefore, the matrix truncated nuclear norm (MTNN) \cite{hu2013fast,ji2017nonconvex} was proposed by minimizing the sum of the  $\min(m,n)-r$ minimum singular values because the rank of a matrix only depends on the first $r$ nonzero singular values. In this way, a more accurate approximation of rank function is obtained, at the same time the empirical research showed that the MTNN approach has much better approximation performances than other methods based matrix nuclear norm \cite{dong2018low}.

Although these low-rank prior based approaches have obtained good results, additional information could be considered for a more accurate reconstruction. Another thing need to be noted is that the low-rank component always indicates that the real data in practice also have intrinsically sparse property \cite{yang2010image,wright2010sparse}. One possible way is to exploit the sparse information of the complete matrix in a certain domain, such as transform domains where many signals have inherent sparse structures \cite{yang2010image}. To describe the sparse property in a certain domain, Dong $\emph{et al.}$  \cite{dong2018low} proposed a general way by applying the transform operation to matrices as an implicit function.

However, dealing with color images and videos by matrix does not exploit the structural information among channels. It is natural to consider extending the matrix completion to tensor completion \cite{long2019lowrank} for such a task. Since there is no perfect definitions for tensor rank and tensor nuclear norm, several types of tensor nuclear norm were proposed. Liu $\emph{et al.}$ \cite{liu2013tensor} initially proposed the sum of matricized nuclear norm (SMNN) of a tensor, which is defined as
\begin{equation}
\begin{aligned}
&\mathop{\rm{min}}\limits_{\bm{\mathcal{X}}}\sum_{i=1}^{n}\|{\bm{\mathcal{X}}_{[i]}}\|_{*} &\\
&s.t. \bm{\mathcal{X}}_{\bm{\Omega}}={\bm{\mathcal{M}}}_{\bm{\Omega}}, &\\
\end{aligned}
\end{equation}
where $\bm{\mathcal{X}}_{[i]}$ denotes the matrix of the tensor unfolded along the $i$th mode, e.g., the mode-$i$ matricization of $\bm{\mathcal{X}}$, $\alpha_{i}>0$ is a parameter which satisfies $\sum_{i=1}^{n}\alpha_{i}=1$, and ${\bm{\mathcal{M}}}_{\bm{\Omega}}$ is the original incomplete tensor. Kilmer $\emph{et al.}$ \cite{kilmer2013third} proposed a novel tensor decomposition method, called the tensor singular value decomposition (T-SVD). Then Zhang $\emph{et al.}$ \cite{zhang2014novel} proposed a new tubal nuclear norm based on T-SVD, which is defined as the sum of nuclear norms of all frontal slices in the Fourier domain and proofed that it was a convex relaxation to the tensor tubal rank. As a result their optimization model can be written as
\begin{equation}
\begin{aligned}
&\mathop{\rm{min}}\limits_{\bm{\mathcal{{X}}}}\sum_{i=1}^{n}\alpha_{i}\|\bm{\bar{X}}^{(i)}\|_{*} &\\
&s.t. \bm{\mathcal{X}}_{\bm{\Omega}}={\bm{\mathcal{M}}}_{\bm{\Omega}}, &\\
\end{aligned}
\end{equation}
where $\bm{\bar{X}}^{(i)}$ will be introduced in the next section.

The same as which has mentioned in matrix completion, the tensor nuclear norm also minimizes all the singular value at the same level which is unfair to the larger singular values, because the larger singular values always contain much more important information. Then tensor truncated nuclear norm was proposed. Han $\emph{et al.}$ \cite{han2017sparse} proposed a tensor truncated nuclear nuclear norm T-TNNS based on MTNN. Xue $\emph{et al.}$ \cite{xue2018lowrank} proposed a tensor truncated nuclear norm T-TNN based on T-SVD, which will be given in next section.

To obtain a more accurate completion performance, we consider the sparse property of the tensor in a transform domain. Here we select the multi-dimensional discrete consine transform (DCT) \cite{zhu2010search}, since signals has a intrinsic sparse property in this transform domain \cite{wright2010sparse}. Further we introduce a $\ell_{1}$-norm regularization term into the objective function to impose local sparsity and to preserve the piecewise smooth property of the reconstructed tensor. Then we solve the objective function by alternating between two steps. The first step is achieved by performing T-SVD to the observed tensor. The second step solves the cost function by the alternating direction method of multipliers (ADMM) \cite{boyd2011distributed}, which is widely used for solving constrained optimization problems because its guarantee of convergence in polynomial time.

The remainder of this paper is organized as follows. Section \ref{sec2} presents the notations and definitions. Section \ref{sec3} gives the proposed new method. Section \ref{sec4} shows the experimental result. Section \ref{sec5} makes a conclusion about this paper.

\section{Notations and Preliminaries}\label{sec2}

In this paper, we denote tensors by boldface Euler script letters, e.g., $\bm{\mathcal{A}}$. Matrices are denoted by boldface capital letters, e.g., $\bm{A}$. Vectors are denoted by boldface lowercase letters, e.g., $\bm{a}$, and scalars are denoted by lowercase letters, e.g., $a$. We denote $\bm{I}_{n}$ as the $n\times n$ identity matrix. The field of real numbers and complex numbers are denoted as $\mathbb{R}$ and $\mathbb{C}$, respectively. For a 3D tensor $\bm{\mathcal{A}}\in\mathbb{C}^{n_1\times n_2\times n_3}$, we denote its $(i,j,k)$-th elements as $\bm{\mathcal{A}}_{ijk}$ or $a_{ijk}$ and use Matlab commands $\bm{\mathcal{A}}(i,:,:)$, $\bm{\mathcal{A}}(:,i,:)$ and $\bm{\mathcal{A}}(:,:,i)$ to respectively denote the $i$-th horizontal, lateral and frontal slice. More often, the frontal slice $\bm{\mathcal{A}}(:,:,i)$ is denoted as $\bm{A}^{(i)}$. The tube is denoted as $\bm{\mathcal{A}}(i,j,:)$. The inner product of $\bm{A}$ and $\bm{B}$ in $\mathbb{R}^{n_1\times n_2}$ is defined as $\langle\bm{A},\bm{B}\rangle=\rm{tr}(\bm{{A}}^{T}\bm{{B}})$, where $\bm{A}^T$ denotes the transpose of $\bm{A}$ and $\rm{tr}(\cdot)$ denotes the matrix trace. The trace of $\bm{\mathcal{A}}$ is defined as ${\rm{tr}}(\bm{\mathcal{A}})=\sum_{i=1}^{n_3}\rm{tr}$$(\bm{{A}}^{(i)})$. The inner product of $\bm{\mathcal{A}}$ and $\bm{\mathcal{B}}$ in $\mathbb{R}^{n_1\times n_2\times n_3}$ is defined as $\langle\bm{\mathcal{A}},{\mathcal{B}}\rangle=\sum_{i=1}^{n_3}\langle\bm{A}^{(i)},\bm{B}^{(i)}\rangle$.

Some norms of tensor and matrix are used. We denote the $\ell_{1}$-norm as $\|\bm{\mathcal{A}}\|_{1}=\sum_{ijk}|a_{ijk}|$ and the Frobenius norm as $\|\bm{\mathcal{A}}\|_{F}= \sqrt{\sum_{ijk}{|a_{ijk}|}^{2}}$. The matrix nuclear norm is $\|\bm{A}\|_{*}=\sum_{i}|\sigma_{i}\bm{(A)}|$, e.g., the sum of all singular values of matrix $\bm{A}$.

For tensor $\bm{\mathcal{A}}\in\mathbb{R}^{n_1\times n_2\times n_3}$, by using the Matlab command fft, we denote $\bm{\bar{\mathcal{A}}}$ as the result of discrete Fourier transform (DFT) \cite{cai2010singular} of $\bm{\mathcal{A}}$ along the third mode, i.e., $\bm{\bar{\mathcal{A}}}= $fft$(\bm{\mathcal{A}},[],3)$. In the same fashion, we can compute $\bm{\mathcal{A}}$ from $\bm{\bar{\mathcal{A}}}$ by ifft$(\bm{\bar{\mathcal{A}}},[],3)$ using the inverse FFT. In particular we denote $\bm{\bar{A}}$ as a block diagonal matrix with each diagonal block as the frontal slice $\bm{\bar{A}}^{(i)}$ of $\bm{\bar{\mathcal{A}}}$, i.e.,
\begin{equation}
\bm{\bar{A}}={\rm{bdiag}}(\bm{\bar{\mathcal{A}}})=
\left[
  \begin{array}{cccc}
  \bm{\bar{A}}^{(1)}& & & \\
  &  \bm{\bar{A}}^{(2)}& & \\
  & & \ddots& \\
  & & & \bm{\bar{A}}^{(n_3)} \\
  \end{array}
  \right].
\end{equation}
This bdiag$(\cdot)$ can be seen as an operator which maps the tensor $\bm{\bar{\mathcal{A}}}$ to the block diagonal matrix $\bm{\bar{A}}$.

The block circulant matrix corresponding to a tensor is defined as
\begin{equation}
\rm{bcirc}(\bm{\mathcal{A}})=
\left[
  \begin{array}{cccc}
  \bm{A}^{(1)}&\bm{A}^{(n_3)} & \ldots&\bm{A}^{(2)} \\
  \bm{A}^{(2)}&\bm{A}^{(1)}& \ldots&\bm{A}^{(3)} \\
  \vdots&\vdots & \ddots&\vdots \\
  \bm{A}^{(n_3)}&\bm{A}^{(n_3-1)} & \ldots& \bm{A}^{(1)} \\
  \end{array}
  \right].
\end{equation}

For tensor $\bm{\mathcal{A}}\in\mathbb{R}^{n_1\times n_2\times n_3}$, we define
\begin{equation}
\rm{unfold}(\bm{\mathcal{A}})=
\left[
  \begin{array}{c}
  \bm{A}^{(1)}\\
  \bm{A}^{(2)}\\
  \vdots \\
  \bm{A}^{(n_3)}\\
  \end{array}
  \right] , {\rm{fold}({\rm{unfold}}(\bm{\mathcal{A}}))}=\bm{\mathcal{A}},
\end{equation}
where the $\rm{unfold}$ operator maps $\bm{\mathcal{A}}$ to a matrix of size $n_1n_3\times n_2$ and $\rm{fold}$ is its inverse operator.

\newtheorem{definition}{\bf Definition}
\newtheorem{thm}{\bf Theorem}

\begin{definition}
{\bf Tensor product} \cite{zhang2014novel} Let $\bm{\mathcal{A}}\in\mathbb{R}^{n_1\times n_2\times n_3}$
and $\bm{\mathcal{B}}\in\mathbb{R}^{n_2\times n_4\times n_3}$. Then the tensor-product $\bm{\mathcal{A}}\ast\bm{\mathcal{B}}$
is defined to be a tensor of size $n_1\times n_4\times n_3$, e.g.,
\begin{equation}
\bm{\mathcal{A}}\ast\bm{\mathcal{B}}=\rm{fold}(\rm{bcirc}(\bm{\mathcal{A}})\cdot \rm{unfold}(\bm{\mathcal{B}})),
\end{equation}
where $\cdot $ denotes the matrix product.
\end{definition}

The tensor product can be understood from two perspectives. First, in the original region, it is analogous to the matrix product except that the
circular convolution replaces the product operation between the elements. The tensor product reduces to the standard matrix product when $n_3=1$.
Second, in Fourier domain, it is equivalent to the matrix multiplication,
e.g., $\bm{\mathcal{C}}=\bm{\mathcal{A}}\ast\bm{\mathcal{B}} \Longleftrightarrow \bm{\bar{C}}=\bm{\bar{A}}\bm{\bar{B}}$ \cite{lu2019tensor} .

\begin{definition}
{\bf Transpose} \cite{zhang2014novel} The transpose of a tensor $\bm{\mathcal{A}}\in\mathbb{R}^{n_1\times n_2\times n_3}$ is
the ${n_2\times n_1\times n_3}$ tensor $\bm{\mathcal{A}}^{T}$ obtained by transposing each of the frontal slice and then reversing
the order of transposed frontal slices 2 through $n_3$, e.g.,
\begin{equation}
\begin{aligned}
&{{(\bm{\mathcal{A}}^{T})}}^{(1)}={(\bm{\mathcal{A}}^{(1)})}^{T} &\\
&{{(\bm{\mathcal{A}}^{T})}}^{(i)}={(\bm{\mathcal{A}}^{(n_3-i+2)})}^{T}, i=2,\ldots,n_3.&\\
\end{aligned}
\end{equation}
\end{definition}

\begin{definition}
{\bf Identity tensor} \cite{zhang2014novel} The identity tensor $\bm{\mathcal{I}}\in\mathbb{R}^{n\times n\times n_3}$ is the tensor whose first frontal slice is the $n\times n$ identity matrix, and whose other frontal slices are all zeros.
\end{definition}

\begin{definition}
{\bf Orthogonal tensor} \cite{zhang2014novel} A tensor $\bm{\mathcal{Q}}\in\mathbb{R}^{n\times n\times n_3}$ is orthogonal if it satisfies
\begin{equation}
{\bm{\mathcal{Q}}}^{T}\ast\bm{\mathcal{Q}}=\bm{\mathcal{Q}}\ast{\bm{\mathcal{Q}}}^{T}=\bm{\mathcal{I}}.
\end{equation}
\end{definition}

\begin{definition}
{\bf F-diagonal tensor} \cite{zhang2014novel} A tensor is called f-diagonal if each of its frontal slices is a diagonal matrix.
\end{definition}

\begin{figure}[htbp]
\centering
\includegraphics[width=9cm]{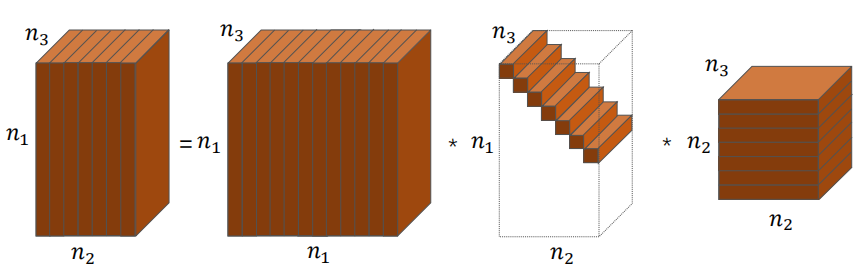}
\caption{Illustration of the T-SVD of an $n_1\times n_2\times n_3$ tensor}\label{fig:tsvd}
\end{figure}

\begin{thm}
{\bf T-SVD} \cite{zhang2014novel} Let $\bm{\mathcal{A}}\in\mathbb{R}^{n_1\times n_2\times n_3}$, then it can be factorized as
\begin{equation}
\bm{\mathcal{A}}=\bm{\mathcal{U}}\ast \bm{\mathcal{S}}\ast {\bm{\mathcal{V}}}^{T},
\end{equation}
where $\bm{\mathcal{U}}\in\mathbb{R}^{n_1\times n_1\times n_3}$, $\bm{\mathcal{V}}\in\mathbb{R}^{n_2\times n_2\times n_3}$ are orthogonal,
and $\bm{\mathcal{S}}\in\mathbb{R}^{n_1\times n_2\times n_3}$ is an f-diagonal tensor. Figure \ref{fig:tsvd} shows an example.
\end{thm}

\begin{definition}\label{E5}
{\bf Tensor tubal rank and tensor nuclear norm} \cite{xue2018lowrank} Let the T-SVD of tensor $\bm{\mathcal{A}}\in\mathbb{R}^{n_1\times n_2\times n_3}$
be $\bm{\mathcal{U}}\ast \bm{\mathcal{S}}\ast {\bm{\mathcal{V}}}^{T}$. The tensor tubal rank of $\bm{\mathcal{A}}$ is defined as the maximum
rank among all frontal slices of the f-diagonal $\bm{\mathcal{S}}$, i.e., max$_{i}$ rank$(\bm{\mathcal{S}}^{(i)})$. The tensor nuclear
norm $\|\bm{\mathcal{A}}\|_{*}$ is defined as the sum of the singular values in all frontal slices of $\bm{\mathcal{S}}$, i.e.,
\begin{equation}\label{E4}
\|{\bm{\mathcal{A}}}\|_{*}={\rm{tr}}(\bm{\mathcal{S}})={\sum_{i=1}^{n_3}}{\rm{tr}}(\bm{S}^{(i)}).
\end{equation}
\end{definition}

According to the definition of FFT function, we carry out fft along the third mode, then we can get a symmetric property between the trace of tensor product
$\bm{\mathcal{A}}\ast \bm{\mathcal{B}}$ and the trace of ${\bm{\bar{A}}}^{(1)}$ and ${\bm{\bar{B}}}^{(1)}$ \cite{xue2018lowrank}, i.e.,
\begin{equation}\label{E3}
\rm{tr}(\bm{\mathcal{A}}\ast\bm{\mathcal{B}})=\rm{tr}(\bm{{\bar{A}}}^{(1)}\bm{{\bar{B}}}^{(1)}).
\end{equation}

According to (\ref{E3}), the tensor nuclear norm defined in (\ref{E4}) can be simplified as \cite{xue2018lowrank}
\begin{equation}\label{E2}
\|\bm{\mathcal{A}}\|_{*}=\rm{tr}(\bm{\mathcal{S}})=\text{tr}(\bm{{\bar{S}}}^{(1)})=\|\bm{{\bar{A}}}^{(1)}\|_{*}.
\end{equation}

The formulation (\ref{E2}) suggests that we can compute the tensor nuclear norm by one matrix SVD in the Fourier domain rather than the complicated T-SVD in the original domain.

\begin{definition}\label{E11}
{\bf Tensor singular value thresholding} \cite{xue2018lowrank} Assume that the T-SVD of tensor $\bm{\mathcal{X}}\in\mathbb{R}^{n_1\times n_2\times n_3}$
is $\bm{\mathcal{U}}\ast\bm{\mathcal{S}}\ast\bm{\mathcal{V}}^{T}$. The singular value thresholding (SVT) \cite{cai2010singular} operator $(\bm{\mathcal{D}}_{\tau})$
is performed on each frontal slice of the f-diagonal tensor $\bm{\mathcal{\bar{S}}}$.
\begin{equation}
\bm{\mathcal{D}}_{\tau}(\bm{\mathcal{X}})=\bm{\mathcal{U}}\ast\bm{\mathcal{D}}_{\tau}(\bm{\mathcal{S}})\ast\bm{\mathcal{V}}^{T},
\end{equation}
\end{definition}
where $\bm{\mathcal{D}}_{\tau}(\bm{\mathcal{S}})$ is the inverse fft of $\bm{\mathcal{D}}_{\tau}(\bm{\mathcal{\bar{S}}})$, and $\bm{\mathcal{D}}_{\tau}(\bm{\mathcal{\bar{S}}}^{(i)}) = {\rm{diag}}({\rm{max}}\{\sigma_{t}-\tau,0\}_{1\leq t\leq r}), i=1,2,\ldots,n_3.$
$\tau > 0$ is a constant, and $r$ is the rank of $(\bm{\bar{S}}^{(i)})$.

\begin{thm}\label{thm1}
 \cite{dong2018low} Let $\bm{X} \in \mathbb{R}^{m\times n}$ be a given matrix and $r$ be any non-negative integer with $r\leq {\rm{min}}(m,n)$. For any
 matrices $\bm{A} \in \mathbb{R}^{r\times m}$, $\bm{B} \in\mathbb{R}^{r\times n}$ satisfying $\bm{AA}^{T}=\bm{I}_{r\times r}$,
 $\bm{BB}^{T}=I_{r\times r}$, we have
\begin{equation}
{\rm{tr}}(\bm{AXB}^{T}) \leq \sum_{i=1}^{r}\sigma_i(\bm{X}),
\end{equation}
where $\bm{I}_{r\times r}$ denotes the identity matrix of size $r\times r$.
\end{thm}

\section{Proposed method }\label{sec3}
In the formulation of our proposed method, the low-rank assumption and the sparse prior are both considered in order to better utilize the structure information of the tensor. Since the truncated nuclear norm can provide a better approximation to rank function in matrix \cite{dong2018low}, Xue \emph{et al}. \cite{xue2018lowrank} extended this property directly to tensor by defining a new tensor truncated nuclear norm (T-TNN). We employ this T-TNN to model the low-rank prior information in this paper. For the sparse prior, we proposed a new term to describe it. For the reasons mentioned in the introduction section,
we assume the original tensor $\bm{\mathcal{X}}$ is sparse in the DCT transform domain. Hence the proposed method is named after sparse regularization in a transformed domain, e.g., SRTD. Let $\bm{\mathcal{T(\cdot)}}$ denote the forward n-dimensianal DCT,
and the transformed tensor
$\bm{\mathcal{E}}=\bm{\mathcal{T(\bm{\mathcal{X}})}}$ is assumed to be sparse.
\subsection{Problem formulation}
 For a tensor $\bm{\mathcal{X}} \in \mathbb{R}^{n_1\times n_2\times n_3}$, the tensor completion problem can be formulated as the following
 constrained optimization problem
\begin{equation}\label{E6}
\begin{aligned}
&\mathop{\rm{min}}\limits_{\bm{\mathcal{X}},\bm{\mathcal{E}}}\|\bm{\mathcal{X}}\|_t+\lambda \|\bm{\mathcal{E}}\|_1 &\\
&s.t. \bm{\mathcal{X}}_{\bm{\Omega}} = \bm{\mathcal{M}}_{\bm{\Omega}} &\\
&\bm{\mathcal{E}} = \bm{\mathcal{T(\mathcal{X})}}, &\\
\end{aligned}
\end{equation}
where $\lambda >0$ and $\bm{\mathcal{M}}_{\bm{\Omega}}$ is the original incomplete tensor with observed values on the support $\bm{\Omega}$. The tensor truncated nuclear norm $\|\bm{\mathcal{X}}\|_t$ can be expressed as follows
 \begin{equation}\label{E17}
 \begin{aligned}
 \|\bm{\mathcal{X}}\|_t &=\|\bar{\bm{X}}^{(1)}\|_t=\sum_{j=r+1}^{min(n_1,n_2)}\sigma_j(\bar{\bm{X}}^{(1)})& \\
 &= \sum_{j=1}^{min(n_1,n_2)}\sigma_j(\bar{\bm{X}}^{(1)})-\sum_{j=1}^{r}\sigma_j(\bar{\bm{X}}^{(1)}). &\\
 \end{aligned}
 \end{equation}
 where $r$ is the truncated singular value.

Since formulation (\ref{E17}) is nonconvex, it is difficult to solve it directly. We use Theorem \ref{thm1} to transform (\ref{E17}) into a convex problem. Combining with (\ref{E3}), (\ref{E2}) and Theorem \ref{thm1}, (\ref{E17}) can be reformulated as
 \begin{equation}
 \begin{aligned}
 \|\bm{\mathcal{X}}\|_t &=\|\bar{\bm{X}}^{(1)}\|_*-\mathop{\rm{max}}\limits_{\begin{array}{c}
 \bar{\bm{A}}^{(1)}\bar{{\bm{A}}}^{{(1)}^{T}}=\bm{I} \\
 \bar{\bm{B}}^{(1)}\bar{{\bm{B}}}^{{(1)}^{T}}=\bm{I}  \\
 \end{array}}
 \rm{tr}(\bar{\bm{{A}}}^{(1)}\bar{\bm{{X}}}^{(1)}\bar{{\bm{{B}}}}^{{(1)}^{T}}) &\\
 &= \|\bm{\mathcal{X}}\|_*-\mathop{\rm{max}}\limits_{\begin{array}{c}
 \bm{\mathcal{A}}\ast\bm{\mathcal{A}}^{T}=\bm{\mathcal{I}} \\
 \bm{\mathcal{B}}\ast\bm{\mathcal{B}}^{T}=\bm{\mathcal{I}} \\
 \end{array}}
 \rm{tr}(\bm{\mathcal{A}}\ast\bm{\mathcal{X}}\ast\bm{\mathcal{B}}^{T}). &\\
 \end{aligned}
 \end{equation}

So the formulation (\ref{E6}) becomes
\begin{equation}\label{E7}
\begin{aligned}
&\mathop{\rm{min}} \limits_{\bm{\mathcal{X}},\bm{\mathcal{E}}}\|\bm{\mathcal{X}}\|_{\ast}-\mathop{\rm{max}}\limits_{\begin{array}{c}
 \bm{\mathcal{A}}\ast\bm{\mathcal{A}}^{T}=\bm{\mathcal{I}} \\
 \bm{\mathcal{B}}\ast\bm{\mathcal{B}}^{T}=\bm{\mathcal{I}} \\
 \end{array}}
 \rm{tr}(\bm{\mathcal{A}}\ast\bm{\mathcal{X}}\ast\bm{\mathcal{B}}^{T})+\lambda \|\bm{\mathcal{E}}\|_1 &\\
&s.t.\bm{\mathcal{X}}_{\bm{\Omega}} = \bm{\mathcal{M}}_{\bm{\Omega}} &\\
&\bm{\mathcal{E}} = \bm{\mathcal{T(\mathcal{X})}}. &\\
\end{aligned}
\end{equation}

To solve the optimization problem (\ref{E7}), an iterative method alternating between two steps is adopted. In the first step, we compute
the T-SVD of a fixed tensor, i.e., $\bm{\mathcal{X}}_k=\bm{\mathcal{U}}\ast\bm{\mathcal{S}}\ast\bm{\mathcal{V}}^{T}$, and then $\bm{\mathcal{A}}_k$
and $\bm{\mathcal{B}}_k$ can be
derived from $\bm{\mathcal{U}}$ and $\bm{\mathcal{V}}$, i.e.,
\begin{equation}
\bm{\mathcal{A}}_k=\bm{\mathcal{U}}(:,1:r,:)^{T},\bm{\mathcal{B}}_k=\bm{\mathcal{V}}(:,1:r,:)^{T}.
\end{equation}
In the second step, by assuming that $\bm{\mathcal{A}}_k$ and $\bm{\mathcal{B}}_k$ are fixed, we compute $\bm{\mathcal{X}}_k$ from a simplified formulation
\begin{equation}\label{E8}
\begin{aligned}
&\mathop{\rm{min}} \limits_{\bm{\mathcal{X}},\bm{\mathcal{E}}}\|\bm{\mathcal{X}}\|_{\ast}-{\rm{tr}}(\bm{\mathcal{A}}_{k}\ast\bm{\mathcal{X}}
\ast\bm{\mathcal{B}}_{k}^{T})
+\lambda \|\bm{\mathcal{E}}\|_1 &\\
&s.t. \bm{\mathcal{X}}_{\bm{\Omega}} = \bm{\mathcal{M}}_{\bm{\Omega}} &\\
&\bm{\mathcal{E}} = \bm{\mathcal{T(\mathcal{X})}}.&\\
\end{aligned}
\end{equation}
ADMM is used to solve (\ref{E8}), and the details will be presented in the next subsection. The overall solution framework for solving (\ref{E7}) is summarized in Algorithm 1.

\begin{algorithm}
\renewcommand{\algorithmicrequire}{\textbf{Input:}}
\renewcommand{\algorithmicensure}{\textbf{Output:}}
\caption{Low rank tensor completion with sparse regularization in a transformed domain.}
\begin{algorithmic}[1]
\REQUIRE $\bm{\mathcal{M}}$, the original incompletion data; $\bm{\Omega}$, the index set of known elements; $\bm{\Omega}^{c}$, the index set of unknown elements; $K$, the maximum iteration number.
\ENSURE ${\bm{\mathcal{X}}}^{k+1}$, the recovered tensor.
\STATE initialize the model parameter, $\bm{\mathcal{X}}^1=\bm{\mathcal{M}}_{\bm{\Omega}}$, $\epsilon = 10^{-3}$, $k=1$, $K=50$;
\STATE repeat until $\|\bm{\mathcal{X}}^{k+1}-\bm{\mathcal{X}}^k\|_{F} \leq\epsilon$ or $k>K$
\STATE Step 1: given $\bm{\mathcal{X}}^k\in\mathbb{R}^{n_1\times n_2\times n_3}$,calculate $$[{\bm{\mathcal{U}}_k,\bm{\mathcal{S}}_k},
\bm{\mathcal{V}}_k]=\text{T-SVD}(\bm{\mathcal{X}}^k),$$ where $\bm{\mathcal{U}}_k\in \bm{\mathbb{R}}^{n_1\times n_1
\times n_3}, \bm{\mathcal{V}}_k\in \bm{\mathbb{R}}^{n_2\times n_2\times n_3}$ are orthogonal tensors.
\STATE Compute $\bm{\mathcal{A}}_k$ and $\bm{\mathcal{B}}_k$ by $$\bm{\mathcal{A}}_k=\bm{\mathcal{U}}_k(:,1:r,:)^{T},\bm{\mathcal{B}}_k=
\bm{\mathcal{V}}_k(:,1:r,:)^{T}.$$
\STATE Step 2: solve the following optimization problem by ADMM
\begin{equation*}
\bm{\mathcal{X}}^{k+1}={\rm{arg}}\mathop{\rm{min}}\limits_{\bm{\mathcal{X}},
\bm{\mathcal{E}}}\|\bm{\mathcal{X}}\|_*-{\rm{tr}}(\bm{\mathcal{A}}_{k}\ast\bm{\mathcal{X}}\ast{\bm{\mathcal{B}}_{k}}^{T})+
\lambda \|\bm{\mathcal{E}}\|_1
\end{equation*}
s.t. $\bm{\mathcal{X}}_{\bm{\Omega}}=\bm{\mathcal{M}}_{\bm{\Omega}}$, and where $\bm{\mathcal{E}}=\bm{\mathcal{T}}(\bm{\mathcal{X}}).$
\STATE $k = k+1.$
\end{algorithmic}
\end{algorithm}

\subsection{Problem reformulation and ADMM}
We introduce an auxiliary tensor $\bm{\mathcal{W}}$ and reformulate the optimization problem (\ref{E8}) as
\begin{equation}\label{E9}
\begin{aligned}
&\mathop{\rm{min}} \limits_{\bm{\mathcal{X}},\bm{\mathcal{E}},\bm{\mathcal{W}}}\|\bm{\mathcal{X}}\|_{\ast}-
{\rm{tr}}(\bm{\mathcal{A}}_{k}\ast\bm{\mathcal{W}}\ast\bm{\mathcal{B}}_{k}^{T})+\lambda \|\bm{\mathcal{E}}\|_1 &\\
&s.t. \bm{\mathcal{X}}_{\bm{\Omega}} = \bm{\mathcal{M}}_{\bm{\Omega}} &\\
&\bm{\mathcal{E}} = \bm{\mathcal{T(\mathcal{X})}}. &\\
&\bm{\mathcal{X}}=\bm{\mathcal{W}} &\\
\end{aligned}
\end{equation}

Due to the introduction of variable $\bm{\mathcal{W}}$, (\ref{E9}) can be addressed by ADMM. The augmented Lagrangian function of (\ref{E9})
becomes
\begin{equation}\label{E10}
\begin{aligned}
\bm{\mathcal{L}}(\bm{\mathcal{X}},\bm{\mathcal{W}},\bm{\mathcal{E}},\bm{\mathcal{Y}},\bm{\mathcal{Z}},\mu)&=  \|\bm{\mathcal{X}}\|_{\ast}-
{\rm{tr}}(\bm{\mathcal{A}}_{k}\ast\bm{\mathcal{W}}\ast
\bm{\mathcal{B}}_{k}^{T})+\lambda \|\bm{\mathcal{E}}\|_{1} &\\
&+\langle\bm{\mathcal{Y}},\bm{\mathcal{X}}-\bm{\mathcal{W}}\rangle+\frac{\mu}{2}\|\bm{\mathcal{X}}-\bm{\mathcal{W}}
\|_{F}^2 &\\
&+\langle\bm{\mathcal{Z}},\bm{\mathcal{E}}-\bm{\mathcal{T}}(\bm{\mathcal{X}})\rangle+\frac{\mu}{2}\|\bm{\mathcal{E}}-
\bm{\mathcal{T}}(\bm{\mathcal{X}})\|_{F}^2, &\\
\end{aligned}
\end{equation}
where $\bm{\mathcal{Y}}$ and $\bm{\mathcal{Z}}$ are Lagrange multiplier tensors of the same size with $\bm{\mathcal{X}}$, and $\mu$ is a
penalty parameter. Based on the basic framework of ADMM, the optimization problem (\ref{E10}) can be solved by alternatively updating one
variable with the others fixed. Specifically, in the $k$th iteration, the variables are updated via the following scheme
\begin{equation}
\left\{
     \begin{array}{lr}
     \bm{\mathcal{X}}^{k+1}={\rm{arg}}\mathop{\rm{min}}\limits_{\bm{\mathcal{X}}}\bm{\mathcal{L}}
     (\bm{\mathcal{X}},\bm{\mathcal{W}}^{k},\bm{\mathcal{E}}^{k},\bm{\mathcal{Y}}^{k},\bm{\mathcal{Z}}^{k},u^k),&\\
     \bm{\mathcal{E}}^{k+1}={\rm{arg}}\mathop{\rm{min}}\limits_{\bm{\mathcal{E}}}\bm{\mathcal{L}}
     (\bm{\mathcal{X}}^{k+1},\bm{\mathcal{W}}^{k},\bm{\mathcal{E}},\bm{\mathcal{Y}}^{k},\bm{\mathcal{Z}}^{k},u^k),&\\
     \bm{\mathcal{W}}^{k+1}={\rm{arg}}\mathop{\rm{min}}\limits_{\bm{\mathcal{W}}}\bm{\mathcal{L}}
     (\bm{\mathcal{X}}^{k+1},\bm{\mathcal{W}},\bm{\mathcal{E}}^{k+1},\bm{\mathcal{Y}}^{k},\bm{\mathcal{Z}}^{k},u^k),&\\
     \bm{\mathcal{Y}}^{k+1}=\bm{\mathcal{Y}}^{k}+\mu^k(\bm{\mathcal{X}}^{k+1}-\bm{\mathcal{W}}^{k+1}),&\\
     \bm{\mathcal{Z}}^{k+1}=\bm{\mathcal{Z}}^{k}+\mu^k(\bm{\mathcal{E}}^{k+1}-\bm{\mathcal{T}}(\bm{\mathcal{X}}^{k+1})),&\\
     \mu^{k+1}={\rm{min}}(\rho \mu^{k},\mu_{{\rm{max}}}).&\\
     \end{array}
     \right.
\end{equation}
where $\rho>1$ is a predetermined constant to increase the penalty, and $\mu_{\rm{max}}$ is a given upper bound for the penalty.

\subsubsection{Update $\bm{\mathcal{X}}^{k+1}$}
\begin{equation}
\begin{aligned}
\bm{\mathcal{X}}^{k+1}&={\rm{arg}}\mathop{\rm{min}}\limits_{\bm{\mathcal{X}}}\bm{\mathcal{L}}
     (\bm{\mathcal{X}},\bm{\mathcal{W}}^{k},\bm{\mathcal{E}}^{k},\bm{\mathcal{Y}}^{k},\bm{\mathcal{Z}}^{k},\mu^{k}) &\\
&={\rm{arg}}\mathop{\rm{min}}\limits_{\bm{\mathcal{X}}}\|\bm{\mathcal{X}}\|_{*}+\langle{\bm{\mathcal{Y}}}^k,\bm{\mathcal{X}}-{\bm{\mathcal{W}}}^k
\rangle
+\frac{\mu^k}{2}\|\bm{\mathcal{X}}-{\bm{\mathcal{W}}}^k\|_{F}^2  &\\
&+\langle{\bm{\mathcal{Z}}}^k,{\bm{\mathcal{E}}}^k-\bm{\mathcal{T}}(\bm{\mathcal{X}})\rangle+\frac{\mu^k}{2}\|{\bm{\mathcal{E}}}^k-
\bm{\mathcal{T}}(\bm{\mathcal{X}})\|_{F}^2 &\\
&={\rm{arg}} \mathop{\rm{min}}\limits_{\bm{\mathcal{X}}}\|\bm{\mathcal{X}}\|_{*}+\frac{\mu^k}{2}\|\bm{\mathcal{X}}-{\bm{\mathcal{W}}}^k+
\frac{{\bm{\mathcal{Y}}}^k}{\mu^k}\|_{F}^2 &\\
&+\frac{\mu^k}{2}\|{\bm{\mathcal{E}}}^k-\bm{\mathcal{T}}(\bm{\mathcal{X}})+\frac{{\bm{\mathcal{Z}}}^k}{\mu^k}\|_{F}^2.\\
\end{aligned}
\end{equation}
Here $\bm{\mathcal{X}}$ cannot be separated from the other variables since the existence of the transform operator $\bm{\mathcal{T}}$ in the
last term. However, the Parseval's theorem \cite{merhav1998approximate} indicates that if the transformation is an unitary under Frobenius norm, the energy of the signal is unchanged.
According to the Parseval's theorem and the unitary invariant property of DCT, the last term can be rewritten as
\begin{equation}
\|{\bm{\mathcal{E}}}^k-\bm{\mathcal{T}}(\bm{\mathcal{X}})+\frac{{\bm{\mathcal{Z}}}^k}{\mu^k}\|_{F}^2 = \|\bm{\mathcal{G}}({\bm{\mathcal{E}}}^k+
\frac{{\bm{\mathcal{Z}}}^k}{\mu^k})-\bm{\mathcal{X}}\|_{F}^2,
\end{equation}
where $\bm{\mathcal{G}}(\cdot)$ denotes the corresponding inverse transform of $\bm{\mathcal{T}}(\cdot)$.

Hence, we have
\begin{equation}
\begin{aligned}
\bm{\mathcal{X}}^{k+1}&={\rm{arg}} \mathop{\rm{min}}\limits_{\bm{\mathcal{X}}}\|\bm{\mathcal{X}}\|_{*}+\frac{\mu^k}{2}\|\bm{\mathcal{X}}-
{\bm{\mathcal{W}}}^k+\frac{{\bm{\mathcal{Y}}}^k}{\mu^k}\|_{F}^2 &\\
&+\frac{\mu^k}{2}\|\bm{\mathcal{G}}({\bm{\mathcal{E}}}^k+\frac{{\bm{\mathcal{Z}}}^k}{\mu^k})-\bm{\mathcal{X}}\|_{F}^2  &\\
&={\rm{arg}} \mathop{\rm{\rm{min}}}\limits_{\bm{\mathcal{X}}}\|\bm{\mathcal{X}}\|_{*} &\\
&+\mu^k\|\bm{\mathcal{X}}-\frac{1}{2}[{\bm{\mathcal{W}}}^k-\frac{{\bm{\mathcal{Y}}}^k}{\mu^k}+\bm{\mathcal{G}}({\bm{\mathcal{E}}}^k
+\frac{{\bm{\mathcal{Z}}}^k}{\mu^k})]\|_{F}^2. &\\
\end{aligned}
\end{equation}
The above problem has a closed-form solution, given by
\begin{equation}\label{E13}
\bm{\mathcal{X}}^{k+1}=\bm{\mathcal{D}}_{\frac{1}{2\mu^{k}}}\{\frac{1}{2}|{{\bm{\mathcal{W}}}^k-\frac{{\bm{\mathcal{Y}}}^k}{\mu^k}+
\bm{\mathcal{G}}({\bm{\mathcal{E}}}^k
+\frac{{\bm{\mathcal{Z}}}^k}{\mu^k})}|\},
\end{equation}
where $\bm{\mathcal{D}}_{\tau}(\cdot)$ is the SVT operator defined in definition \ref{E11}.

\subsubsection{Update $\bm{\mathcal{E}}^{k+1}$}
\begin{equation}
\begin{aligned}
\bm{\mathcal{E}}^{k+1}&={\rm{arg}}\mathop{\rm{min}}\limits_{\bm{\mathcal{E}}}\mathcal{L}
     (\bm{\mathcal{X}}^{k+1},\bm{\mathcal{W}}^{k},\bm{\bm{\mathcal{E}}},\bm{\mathcal{Y}}^{k},\bm{\mathcal{Z}}^{k},\mu^k) &\\
&={\rm{arg}}\mathop{\rm{min}}\limits_{\bm{\mathcal{E}}}\lambda\|\bm{\mathcal{E}}\|_{1}+\langle{\bm{\mathcal{Z}}}^k,\bm{\mathcal{E}}-
\bm{\mathcal{T}}({\bm{\mathcal{X}}}^{k+1})\rangle &\\
&+\frac{\mu^k}{2}\|\bm{\mathcal{E}}-\bm{\mathcal{T}}({\bm{\mathcal{X}}}^{k+1})\|_{F}^2 &\\
&={\rm{arg}}\mathop{\rm{min}}\limits_{\bm{\mathcal{E}}}\lambda\|\bm{\mathcal{E}}\|_{1}+\frac{\mu^k}{2}\|\bm{\mathcal{E}}-
\bm{\mathcal{T}}({\bm{\mathcal{X}}}^{k+1})+\frac{{\bm{\mathcal{Z}}}^k}{\mu^k}\|_{F}^2. &\\
\end{aligned}
\end{equation}
The above problem has a closed-form solution, given by
\begin{equation}\label{E14}
\bm{\mathcal{E}}^{k+1}=\bm{\mathcal{S}}_{\frac{\lambda}{\mu^k}}(\bm{\mathcal{T}}({\bm{\mathcal{X}}}^{k+1})-\frac{{\bm{\mathcal{Z}}}^k}{\mu^k}),
\end{equation}
where $\bm{\mathcal{S}}_{\tau}(\cdot)$ is the element-wise soft thresholding operator \cite{cai2010singular}, defined by
\begin{equation}
\bm{\mathcal{S}}_{\tau}(x)={\rm}{sgn}(x)\cdot {\rm}{max}\{|x|-\tau,0\}.
\end{equation}

\subsubsection{Update $\bm{\mathcal{W}}^{k+1}$}
\begin{equation}\label{E12}
\begin{aligned}
\bm{\mathcal{W}}^{k+1}&={\rm{arg}}\mathop{\rm{min}}\limits_{\bm{\mathcal{W}}}\bm{\mathcal{L}}
     (\bm{\mathcal{X}}^{k+1},\bm{\mathcal{W}},\bm{\mathcal{E}}^{k+1},\bm{\mathcal{Y}}^{k},\bm{\mathcal{Z}}^{k},\mu^k) & \\
&={\rm{arg}}\mathop{\rm{\rm{min}}}\limits_{\bm{\mathcal{W}}}-{\rm{tr}}(\bm{\mathcal{A}}_{k}\ast\bm{\mathcal{W}}\ast\bm{\mathcal{B}}_{k}^{T})
& \\
&+\langle{\bm{\mathcal{Y}}}^k,{\bm{\mathcal{X}}}^{k+1}-\bm{\mathcal{W}}\rangle +\frac{\mu^k}{2}\|{\bm{\mathcal{X}}}^{k+1}-\bm{\mathcal{W}}\|_{F}^2. &\\
\end{aligned}
\end{equation}
Therefore, by setting the derivative of (\ref{E12}) to zero, we obtain a closed-form solution as follows:
\begin{equation}\label{E15}
\begin{aligned}
&-\bm{\mathcal{A}}_{k}^{T}\ast\bm{\mathcal{B}}_{k}-\bm{\mathcal{Y}}-\mu^{k}({\bm{\mathcal{X}}}^{k+1}-\bm{\mathcal{W}})=0,  &\\
&\bm{\mathcal{W}}^{k+1}=\bm{\mathcal{X}}^{k+1}+\frac{1}{\mu^k}({\bm{\mathcal{A}}_{k}}^T\ast\bm{\mathcal{B}}_{k}+\bm{\mathcal{Y}}^{k}). &\\
\end{aligned}
\end{equation}
In addition, the observed data should keep constant in each iteration, i.e.,
\begin{equation}\label{E16}
\bm{\mathcal{W}}^{k+1}=\bm{\mathcal{X}}^{k+1}_{{\Omega}^c}+\bm{\mathcal{M}}_{\Omega}.
\end{equation}

\subsubsection{Update $\bm{\mathcal{Y}}^{k+1}$}
\begin{equation}\label{E18}
\bm{\mathcal{Y}}^{k+1}=\bm{\mathcal{Y}}^{k}+\mu^k({\bm{\mathcal{X}}^{k+1}-\bm{\mathcal{W}}^{k+1}}).
\end{equation}
\subsubsection{Update $\bm{\mathcal{Z}}^{k+1}$}
\begin{equation}\label{E19}
\bm{\mathcal{Z}}^{k+1}=\bm{\mathcal{Z}}^{k}+\mu^k({\bm{\mathcal{E}}^{k+1}-\bm{\mathcal{T}}(\bm{\mathcal{X}}^{k+1})}).
\end{equation}
\subsubsection{Update $\mu^{k+1}$}
\begin{equation}\label{E20}
\mu^{k+1}={\rm{min}}(\rho \mu^{k},\mu_{{\rm{max}}}).
\end{equation}

The whole procedure to solve the problem (\ref{E9}) is summarized in Algorithm 2.

\begin{algorithm}
\renewcommand{\algorithmicrequire}{\textbf{Input:}}
\renewcommand{\algorithmicensure}{\textbf{Output:}}
\caption{The optimization algorithm to solve the problem (\ref{E9}) by ADMM.}
\begin{algorithmic}[1]
\REQUIRE $\bm{\mathcal{M}}$, the original incompletion data; $\bm{\Omega}$, the index set of known elements; $\epsilon$, a small positive threshold; $\lambda$, a positive hyperparameter; $\mu_{\rm{max}}$, maximum penalty; $K$, maximum iteration number.
\ENSURE the recovered tensor $\bm{\mathcal{X}}.$

\STATE Initialize the model parameters, $k=1$, $\bm{\mathcal{X}}^k=\bm{\mathcal{W}}^{k}=\bm{\mathcal{M}}$,$\bm{\mathcal{E}=0}$, let $\bm{\mathcal{Y}}^{k}$ be a random tensor with the size same as $\bm{\mathcal{X}}^k$, and let $\mu^1$ be a small initial penalty.
\STATE Update $\bm{\mathcal{X}}^{k+1}$ by equations (\ref{E13}),
\STATE Update $\bm{\mathcal{E}}^{k+1}$ by equations (\ref{E14}),
\STATE Update $\bm{\mathcal{W}}^{k+1}$ by equations (\ref{E15}) and (\ref{E16}),
\STATE Update $\bm{\mathcal{Y}}^{k+1}$ by equations (\ref{E18}),
\STATE Update $\bm{\mathcal{Z}}^{k+1}$ by equations (\ref{E19}),
\STATE Update $\mu^{k+1}$ by equations (\ref{E20}),
\STATE If $\|\bm{\mathcal{X}}^{k+1}-\bm{\mathcal{X}}^k\|_{F} \leq\epsilon$, let $\bm{\mathcal{X}}=\bm{\mathcal{X}}^{k+1}$ and stop the iteration.
Otherwise set $k=k+1$ and return to step 2.

\end{algorithmic}
\end{algorithm}
\section{Experiments}\label{sec4}
In this section, several experiments are conducted to demonstrate the efficiency of proposed SRTD method. The compared methods are:
\begin{enumerate}
\item Matrix completion by MTNN \cite{dong2018low};
\item Tensor completion by T-TNN \cite{xue2018lowrank};
\item Tensor completion by T-TNNS \cite{han2017sparse};
\item Tensor completion by SRTD [Ours];
\end{enumerate}

It is necessary to explain the difference between proposed SRTD method and compared methods: MTNN considers transforming the tensor data to matrix data during experiments, which doesn't employ the correlation between channels; T-TNN employs tensor truncated nuclear norm defined by T-SVD, which only considers the low rank information of recovered data; T-TNNS considers the tensor truncated nuclear norm which is defined by the sum of matricized nuclear norm and sparse regularization, but it may can't utilize the channels information as well as the proposed SRTD.

All the experiments are performed in Matlab R2016a on Windows 10, with an Intel Core i5 CPU @2.50GHz and 8 GB Memory.

The Peak Signal-to-Noise  ration (PSNR) is used to describe and evaluate the performance of the recovered images and videos,
which is defined as follows
\begin{equation}
\rm{MSE}=\frac{\|(\bm{\mathcal{X}}_{rec}-\bm{\mathcal{M}})_{\Omega^{c}}\|_{F}}{T},
\end{equation}
\begin{equation}
\rm{PSNR}=10\times log_{10}(\frac{255^2}{MSE})(dB),
\end{equation}
where $\rm{T}$ is the total number of data in a tensor, i.e. $T=n_1n_2n_3$ for tensors considered in Definition 1, and we assume that the maximum pixel value in $\bm{\mathcal{X}}$ is 255. It is obviously that the higher the PSNR, the better the recovery performance.

\subsection{Parameter setting}
\begin{figure}[htbp]
\centering
\includegraphics[width=2.4in]{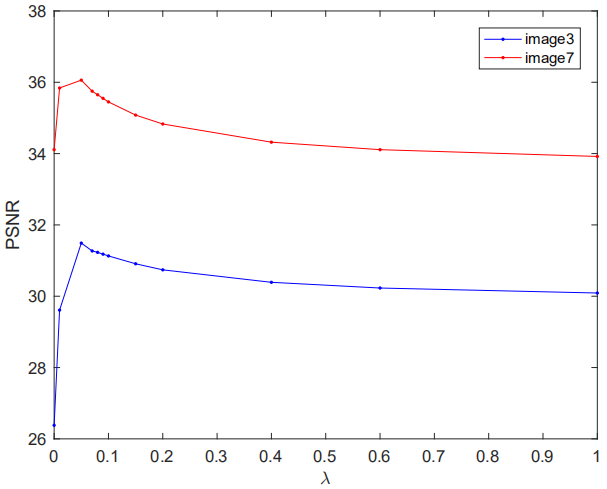}
\caption{PSNR with different $\lambda$.}\label{Fig:2}
\end{figure}

To make sure the comparison is fair, we choose the best parameter for each algorithm. For MTNN,
the parameters are set as $\lambda=0.1$, $\beta=10^{-3}$, $r=15$, and $\epsilon=10^{-3}$, which have been discussed in \cite{dong2018low}.
 For T-TNN, to better illustrate its performance, the random sampling rata (SR) at 50\% is tested. For T-TNNS, the parameter $\lambda$ is set as 0.19 to obtain it's best performance, and the other parameters are set as the same with \cite{han2017sparse}. For the proposed SRTD, there is another parameter $\lambda$ need to be discussed. We test for
 $$\lambda=[0,0.01,0.05,0.07,0.08,0.09,0.1,0.15,0.2,0.4,0.6,1]$$ at SR = 50\%.
 And the PSNR of image 3 and image 7 in Table \ref{table1} are shown in Figure \ref{Fig:2}.

From Figure \ref{Fig:2}, we can see that when $\lambda = 0$, the optimization function of SRTD is reduced to the objection function of T-TNN. We can directly see that the PSNR of SRTD method is better than T-TNN method when $\lambda$ range into (0,1). The PSNR reaches its peak when $\lambda$ is around 0.05, so we set $\lambda=0.05$ in our tests. We do not know the real rank of the incomplete tensor and there is no prior information for us to determine the number of truncated singular values, so following a common practice we manually test the rank range into (1,20) to find the best value in each case.

\subsection{Image recovery with random mask}
A color image can be seen as a 3D tensor usually with a low-rank structure, we first consider ten color images with size of $300\times 400\times 3$. The test sampling rates (SRs) are set as 30\%, 40\% and 50\%. We show the completion result in the Table \ref{table1}. To verify the efficiency of proposed SRTD method, we further randomly select ten color images from Berkeley Segmentation database\footnote{http://www.eecs.berkeley.edu/Research/Projects/CS/vision/bsds/} with size $321\times 481\times 3$. We test the same SRs at 30\%, 40\% and 50\%, and show the completed results in the Table \ref{table2}.

\begin{table}
\centering
\caption{The PSNR value (dB) of the first ten images}
\begin{tabular}{|c|c|c|c|c|c|c|}
\hline
\multirow{2}{*}{No.} &\multirow{2}{*}{Images} & \multirow{2}{*}{SR} & \multicolumn{4}{c|}{PSNR}\\
\cline{4-7}
& & &MTNN &T-TNN &T-TNNS &SRTD\\
\hline
\multirow{3}{*}{1}&\multirow{3}{*}{\includegraphics[width=1.6cm]{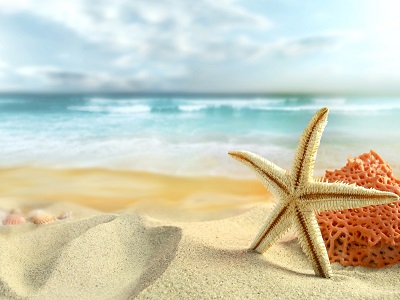}}&30\% &26.48 & 26.35 &27.49 &$\bm{28.85}$ \\
& &40\% &26.86 & 27.96 &28.82 &$\bm{30.37}$\\
& &50\% &27.12 & 29.56 &30.09 &$\bm{31.77}$\\
\hline
\multirow{3}{*}{2}&\multirow{3}{*}{\includegraphics[width=1.6cm]{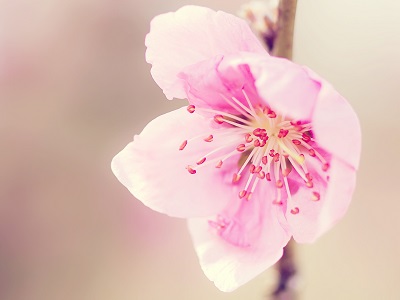}}&30\% &30.89 & 29.22 &31.28 &$\bm{32.58}$ \\
& &40\% &31.68 & 30.80 &32.75 &$\bm{34.11}$\\
& &50\% &32.19 & 32.32 &34.11 &$\bm{35.47}$\\

\hline
\multirow{3}{*}{3}&\multirow{3}{*}{\includegraphics[width=1.6cm]{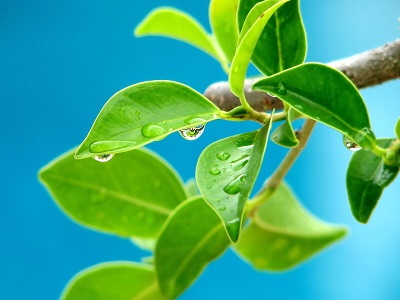}}&30\% &26.31 & 23.32 &26.95 &$\bm{27.76}$ \\
& &40\% &27.21& 25.58 &28.82 &$\bm{29.71}$\\
& &50\% &28.06 & 27.62 &30.43 &$\bm{31.49}$\\

\hline
\multirow{3}{*}{4}&\multirow{3}{*}{\includegraphics[width=1.6cm]{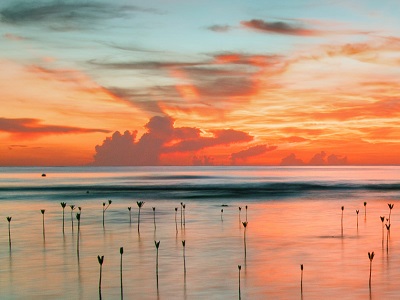}}&30\% &30.27 & 30.74 &31.10 &$\bm{33.08}$ \\
& &40\% &31.20 & 33.02 &33.18 &$\bm{35.17}$\\
& &50\% &32.01 & 35.11 &35.05 &$\bm{36.93}$\\
\hline

\multirow{3}{*}{5}&\multirow{3}{*}{\includegraphics[width=1.6cm]{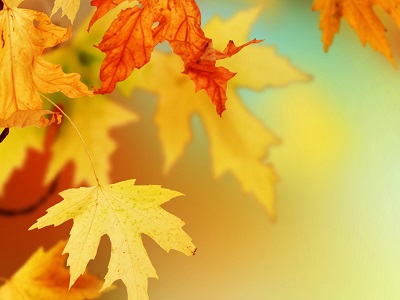}}&30\% &29.98 & 27.77 &30.30 &$\bm{31.56}$ \\
& &40\% &30.67 & 29.50 &31.70 &$\bm{32.78}$\\
& &50\% &31.22& 31.34 &33.07 &$\bm{34.30}$\\
\hline

\multirow{3}{*}{6}&\multirow{3}{*}{\includegraphics[width=1.6cm]{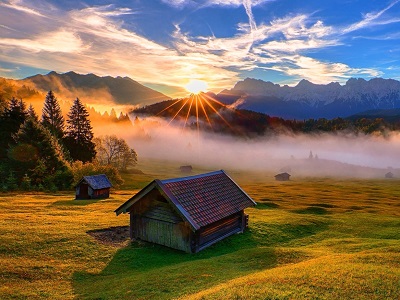}}&30\% &24.02 & 23.28 &24.88 &$\bm{25.95}$ \\
& &40\% &24.35 & 24.60 &26.04 &$\bm{27.20}$\\
& &50\% &24.63 & 25.95 &27.25 &$\bm{28.57}$\\
\hline

\multirow{3}{*}{7}&\multirow{3}{*}{\includegraphics[width=1.6cm]{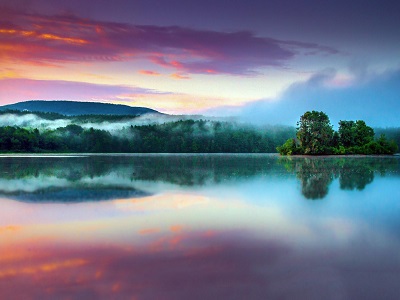}}&30\% &30.64 & 30.75 &31.56 &$\bm{33.14}$\\
& &40\% &31.23 & 32.68 &32.91 &$\bm{34.59}$\\
& &50\% &31.79 & 34.42 &34.37 &$\bm{36.06}$\\
\hline

\multirow{3}{*}{8}&\multirow{3}{*}{\includegraphics[width=1.6cm]{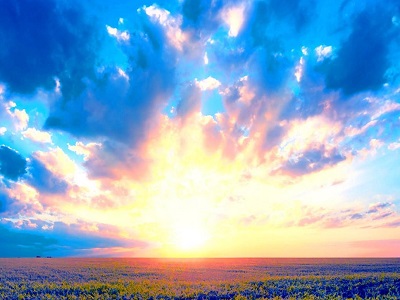}}&30\% &28.05 & 26.39 &28.18 &$\bm{29.03}$ \\
& &40\% &28.59 & 28.11 &29.53 &$\bm{30.55}$\\
& &50\% &29.07 & 29.83 &30.85 &$\bm{32.10}$\\
\hline

\multirow{3}{*}{9}&\multirow{3}{*}{\includegraphics[width=1.6cm]{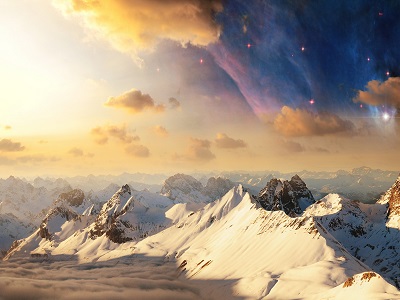}}&30\% &25.15 & 25.21 &26.45 &$\bm{27.75}$ \\
& &40\% &25.46 & 26.90 &27.91 &$\bm{29.45}$\\
& &50\% &25.85 & 28.81 &29.56 &$\bm{31.25}$\\
\hline

\multirow{3}{*}{10}&\multirow{3}{*}{\includegraphics[width=1.6cm]{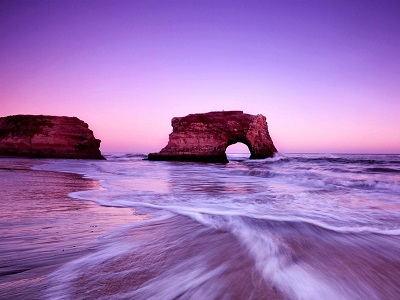}}&30\% &27.43 & 27.97 &29.41&$\bm{31.31}$\\
& &40\% &28.11 & 30.04 &31.34 &$\bm{33.18}$\\
& &50\% &28.74 & 32.29 &33.33 &$\bm{35.23}$\\
\hline
\end{tabular}\label{table1}
\end{table}

\begin{table}
\centering
\caption{The PSNR value (dB) of the ten images from Berkeley Segmentation database}
\begin{tabular}{|c|c|c|c|c|c|c|}
\hline
\multirow{2}{*}{No.} &\multirow{2}{*}{Images} & \multirow{2}{*}{SR} & \multicolumn{4}{c|}{PSNR}\\
\cline{4-7}
& & &MTNN &T-TNN &T-TNNS &SRTD\\
\hline
\multirow{3}{*}{1}&\multirow{3}{*}{\includegraphics[width=1.7cm]{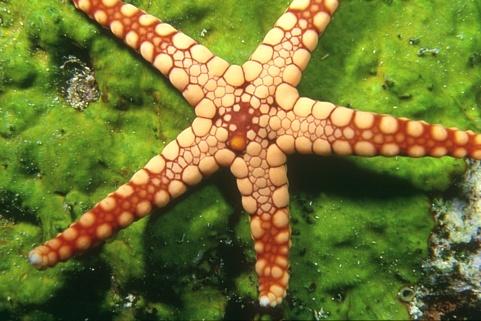}}&30\% &24.35 & 24.00 &26.75 &$\bm{28.24}$ \\
& &40\% &25.40 & 26.64 &29.16 &$\bm{30.81}$\\
& &50\% &26.38 & 29.38 &31.50 &$\bm{33.26}$\\
\hline
\multirow{3}{*}{2}&\multirow{3}{*}{\includegraphics[width=1.7cm]{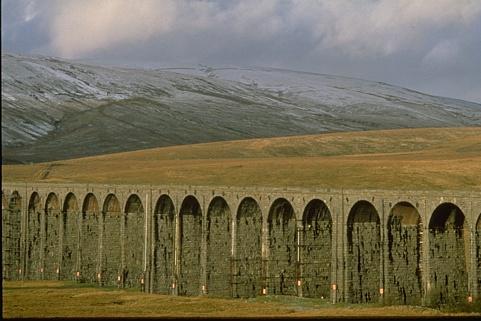}}&30\% &26.43 & 28.68 &28.24 &$\bm{30.44}$ \\
& &40\% &26.98 & 30.87 &30.37 &$\bm{32.70}$\\
& &50\% &27.53 & 33.48 &32.80 &$\bm{35.07}$\\

\hline
\multirow{3}{*}{3}&\multirow{3}{*}{\includegraphics[width=1.7cm]{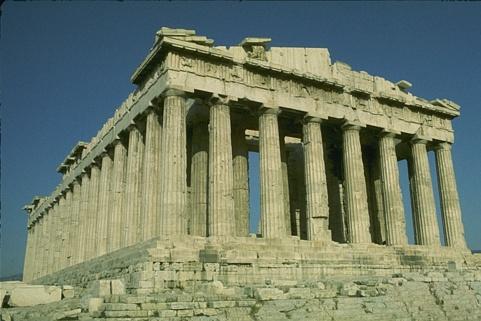}}&30\% &24.05 & 25.62 &26.37 &$\bm{28.52}$ \\
& &40\% &24.82 & 28.50 &28.97 &$\bm{31.21}$\\
& &50\% &25.56 & 30.89 &31.47 &$\bm{33.72}$\\

\hline
\multirow{3}{*}{4}&\multirow{3}{*}{\includegraphics[width=1.7cm,height=1cm]{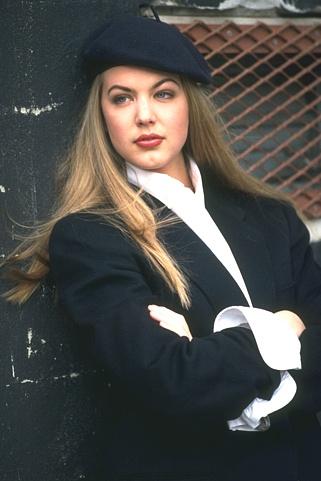}}&30\% &25.60 & 26.67 &29.33&$\bm{30.93}$ \\
& &40\% &26.95 & 26.36 &31.41 &$\bm{33.22}$\\
& &50\% &27.82 & 32.54 &33.87 &$\bm{35.87}$\\
\hline

\multirow{3}{*}{5}&\multirow{3}{*}{\includegraphics[width=1.7cm]{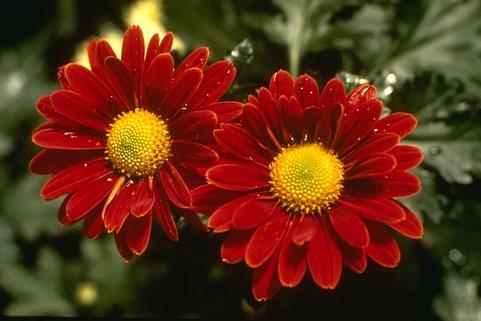}}&30\% &27.07 & 25.97 &28.53 &$\bm{29.78}$ \\
& &40\% &27.97 & 28.37 &30.58 &$\bm{32.02}$\\
& &50\% &28.77 & 30.62 &32.35 &$\bm{33.83}$\\
\hline

\multirow{3}{*}{6}&\multirow{3}{*}{\includegraphics[width=1.7cm]{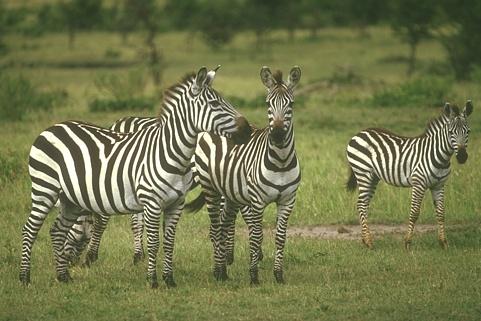}}&30\% &20.11 & 23.01 &22.58 &$\bm{24.89}$ \\
& &40\% &20.85& 25.51 &24.99 &$\bm{27.46}$\\
& &50\% &21.46 & 28.53 &27.82 &$\bm{30.15}$\\
\hline

\multirow{3}{*}{7}&\multirow{3}{*}{\includegraphics[width=1.7cm]{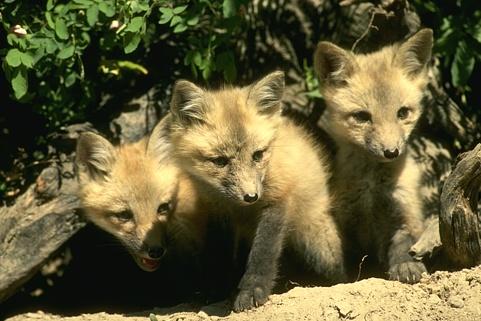}}&30\% &23.47 & 24.46 &26.54 &$\bm{28.24}$\\
& &40\% &24.45 & 27.07 &28.82 &$\bm{30.59}$\\
& &50\% &25.33 & 29.76 &31.05 &$\bm{32.89}$\\
\hline

\multirow{3}{*}{8}&\multirow{3}{*}{\includegraphics[width=1.7cm,height=1cm]{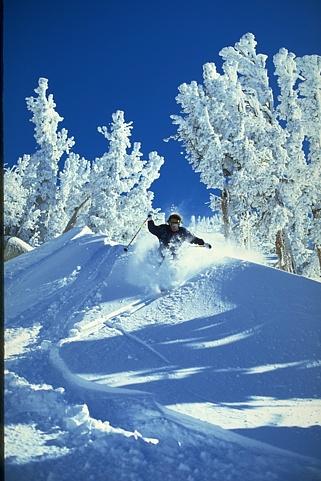}}&30\% &20.56 & 21.56 &22.30 &$\bm{23.93}$\\
& &40\% &21.01 & 23.51 &24.10 &$\bm{25.97}$\\
& &50\% &21.44 & 25.63 &26.21 &$\bm{28.17}$\\
\hline

\multirow{3}{*}{9}&\multirow{3}{*}{\includegraphics[width=1.7cm]{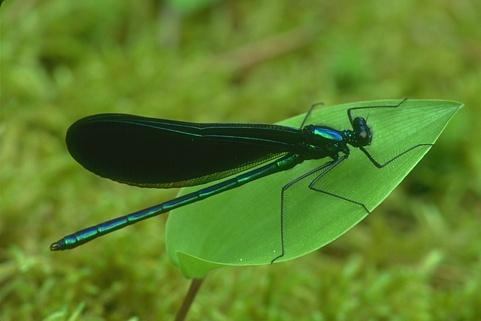}}&30\% &30.65 & 29.53 &32.88 &$\bm{34.21}$ \\
& &40\% &31.45 & 31.89 &34.93 &$\bm{36.41}$\\
& &50\% &32.31 & 34.54 &36.94 &$\bm{38.51}$\\
\hline

\multirow{3}{*}{10}&\multirow{3}{*}{\includegraphics[width=1.7cm]{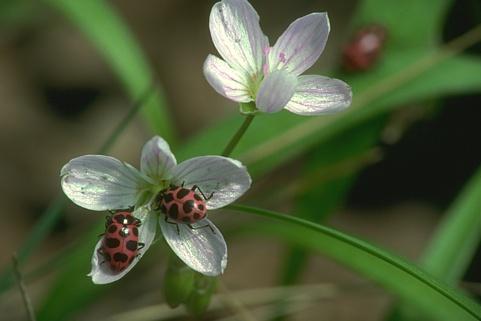}}&30\% &29.88 & 29.92 &32.32 &$\bm{33.92}$ \\
& &40\% &30.78 &32.59 &34.49 &$\bm{36.31}$\\
& &50\% &31.66 & 35.36 &36.69&$\bm{38.65}$\\
\hline
\end{tabular}\label{table2}
\end{table}

 We can easily see the proposed SRTD performances better than compared methods from Table \ref{table1} and Table \ref{table2}. Moreover, in most cases the PSNR of T-TNN, T-TNNS are better than MTNN, which shows that it is better to use tensors rather than matrices to deal with color images.

\begin{figure*}
\centering
\subfigure{
\centering
\includegraphics[width=1in]{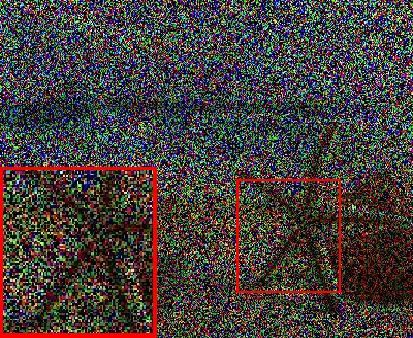}
}
\subfigure{
\centering
\includegraphics[width=1in]{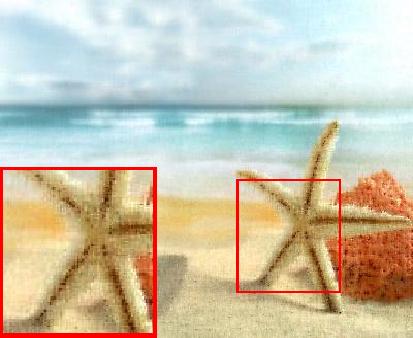}
}
\subfigure{
\includegraphics[width=1in]{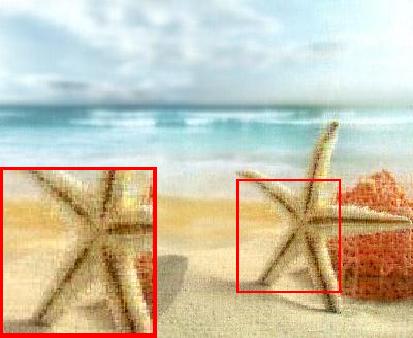}
}
\subfigure{
\includegraphics[width=1in]{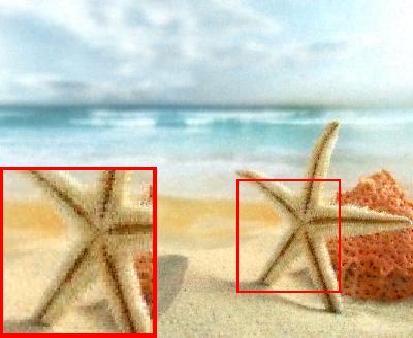}
}
\subfigure{
\includegraphics[width=1in]{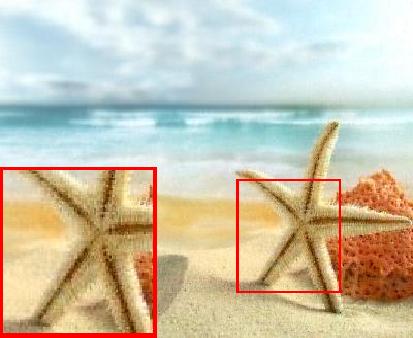}
}
\quad
\subfigure{
\centering
\includegraphics[width=1in]{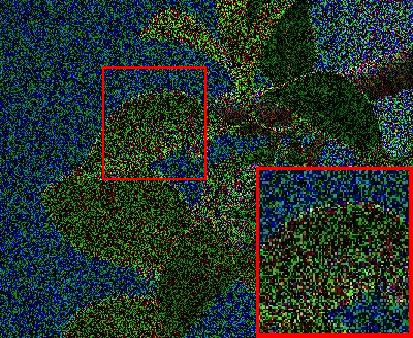}
}
\subfigure{
\centering
\includegraphics[width=1in]{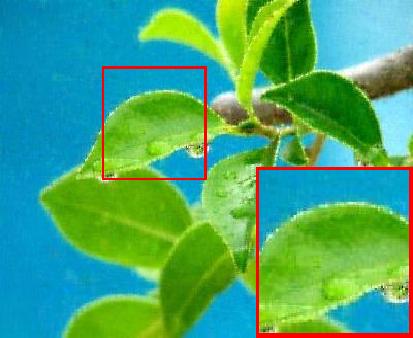}
}
\subfigure{
\includegraphics[width=1in]{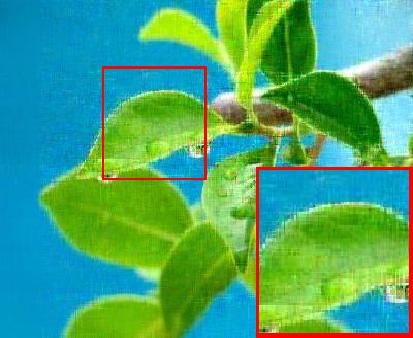}
}
\subfigure{
\includegraphics[width=1in]{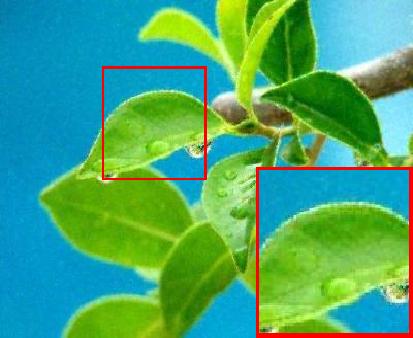}
}
\subfigure{
\includegraphics[width=1in]{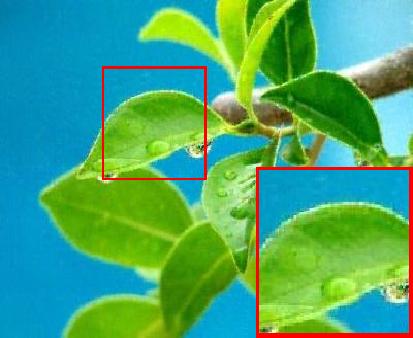}
}
\quad
\subfigure{
\centering
\includegraphics[width=1in]{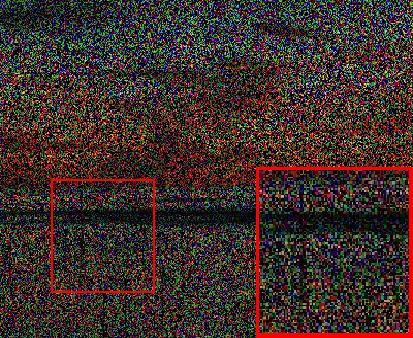}
}
\subfigure{
\centering
\includegraphics[width=1in]{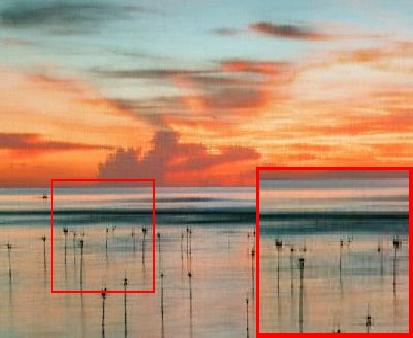}
}
\subfigure{
\includegraphics[width=1in]{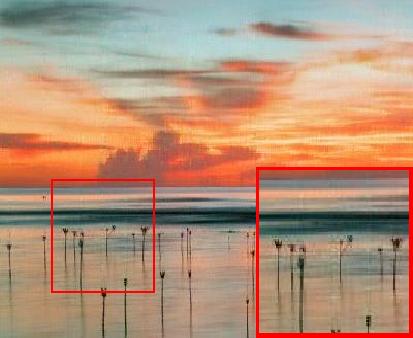}
}
\subfigure{
\includegraphics[width=1in]{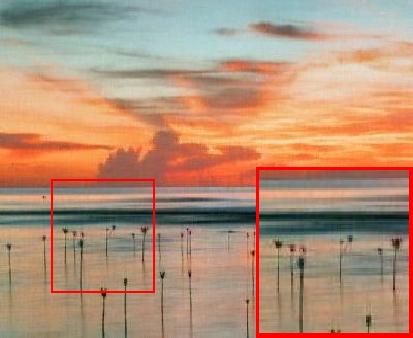}
}
\subfigure{
\includegraphics[width=1in]{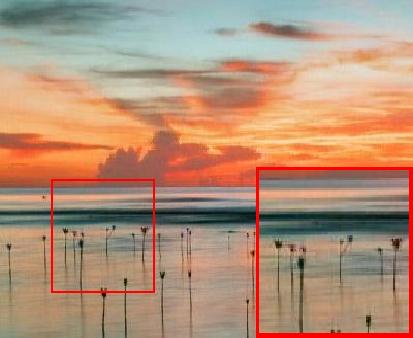}
}
\quad
\subfigure{
\centering
\includegraphics[width=1in]{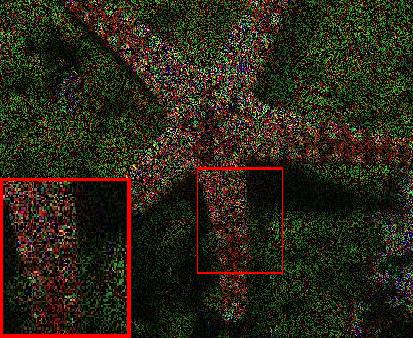}
}
\subfigure{
\centering
\includegraphics[width=1in]{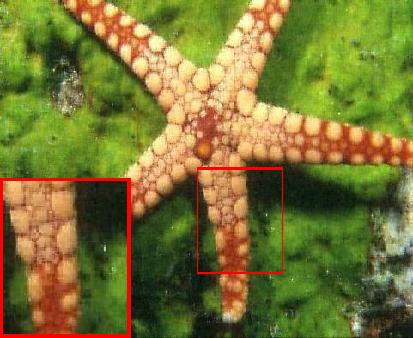}
}
\subfigure{
\includegraphics[width=1in]{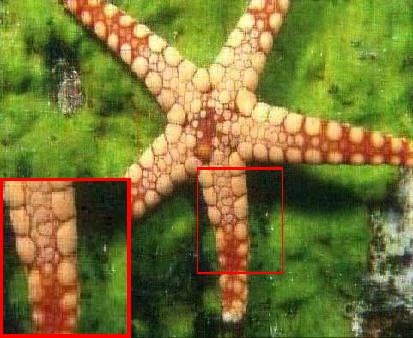}
}
\subfigure{
\includegraphics[width=1in]{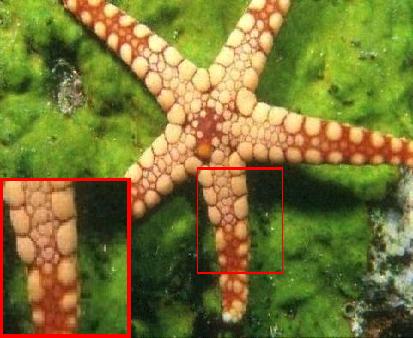}
}
\subfigure{
\includegraphics[width=1in]{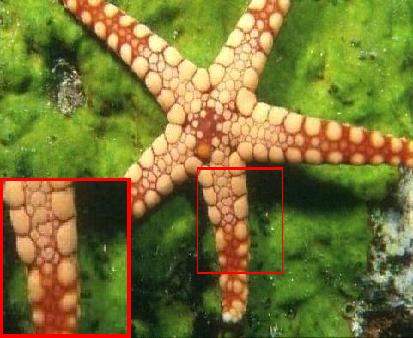}
}
\centering
\subfigure{
\centering
\includegraphics[width=1in]{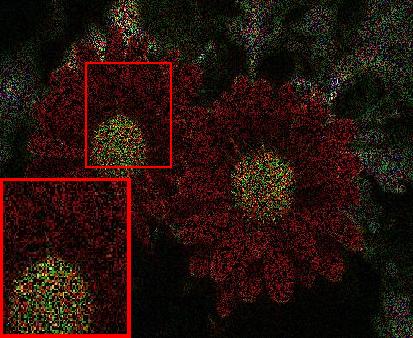}
}
\subfigure{
\centering
\includegraphics[width=1in]{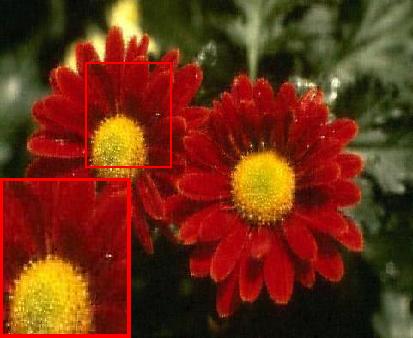}
}
\subfigure{
\includegraphics[width=1in]{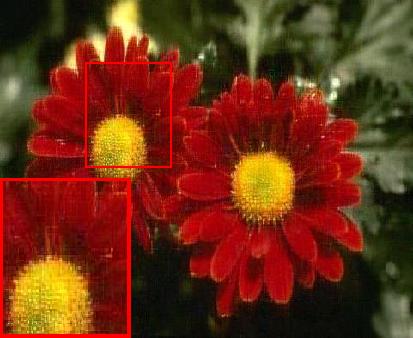}
}
\subfigure{
\includegraphics[width=1in]{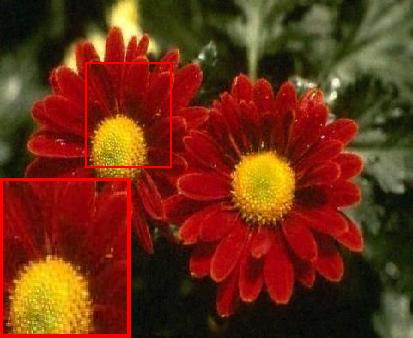}
}
\subfigure{
\includegraphics[width=1in]{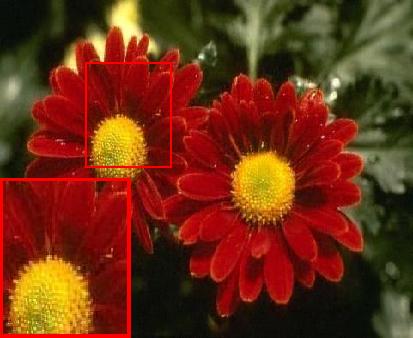}
}
\centering
\subfigure{
\centering
\includegraphics[width=1in]{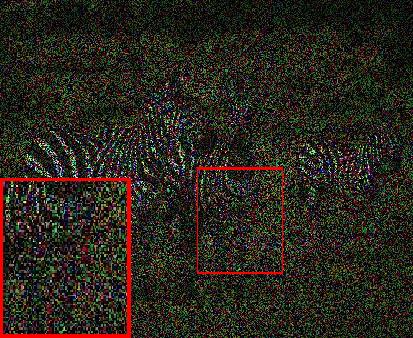}
}
\subfigure{
\centering
\includegraphics[width=1in]{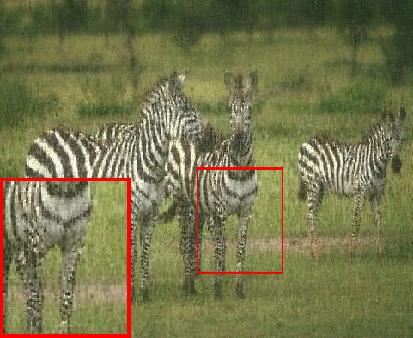}
}
\subfigure{
\includegraphics[width=1in]{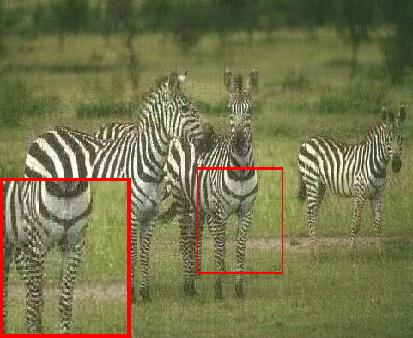}
}
\subfigure{
\includegraphics[width=1in]{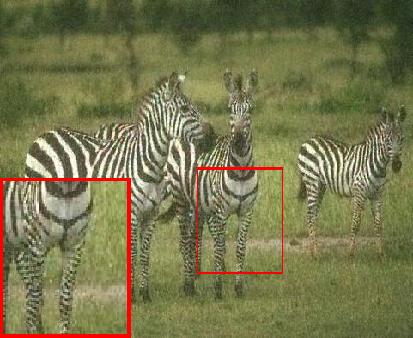}
}
\subfigure{
\includegraphics[width=1in]{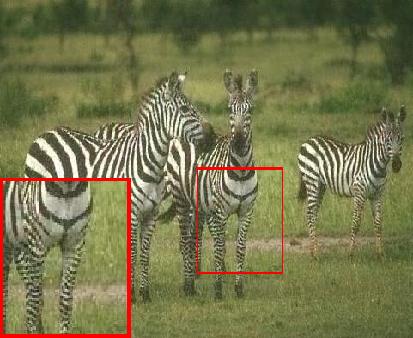}
}
\caption{The result recovered by MTNN, T-TNN, T-TNNS and SRTD at SR = 30\%.}\label{Fig:4}
\end{figure*}

 Figure \ref{Fig:4} shows the six test color images recovered by MTNN, T-TNN, T-TNNS and SRTD respectively at SR $=30\%$. To better show the details of the images, we magnify a significant region for each completed images. Both from the visual quality and the value of PSNR, we can see that the proposed SRTD has a better performance.

\begin{figure}[htbp]
\centering
\subfigure[Image 3]{
\centering
\includegraphics[width=2.3in]{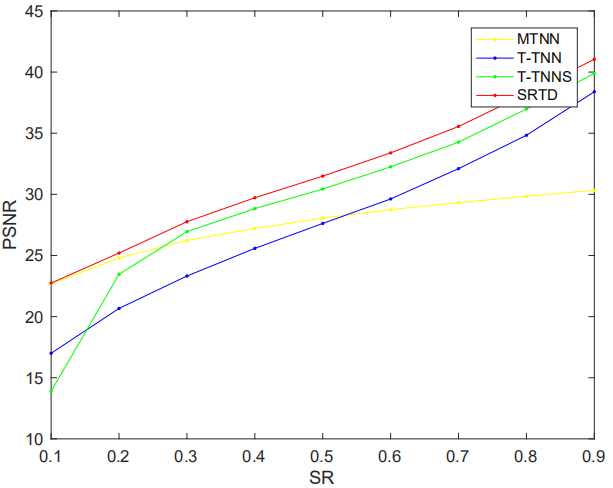}
}
\quad
\subfigure[Image 4]{
\centering
\includegraphics[width=2.3in]{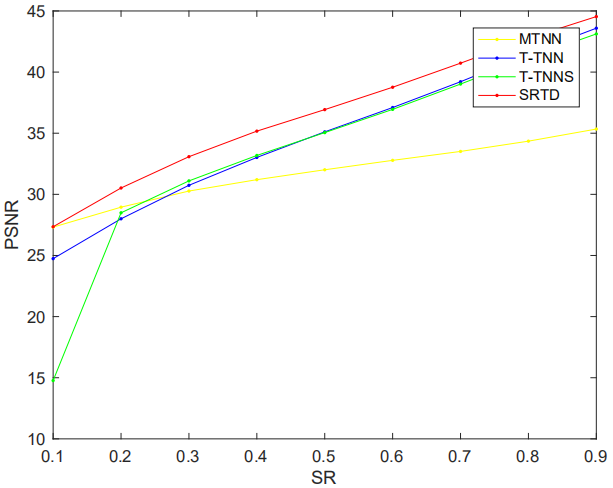}
}
\centering
\caption{The PSNR of image 3 and image 4 for different SR.}\label{Fig:6}
\end{figure}

We further choose image 3 and image 4 from Table \ref{table1}, and we set the SR range from 10\% to 90\%. The result is shown in Figure \ref{Fig:6},
we can directly see from the PSNR value that the performance of SRTD is always better than that of MTNN, T-TNN and T-TNNS, even when the SR becomes very low. However, with the SR decreasing, the performance of T-TNNS decays quickly.

\subsection{Image recovery with text mask}
\begin{figure*}[htbp]
\centering
\subfigure[Original image.]{
\centering
\includegraphics[width=1.5in]{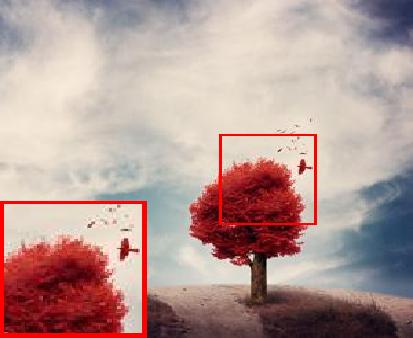}
}
\subfigure[text mask image.]{
\centering
\includegraphics[width=1.5in]{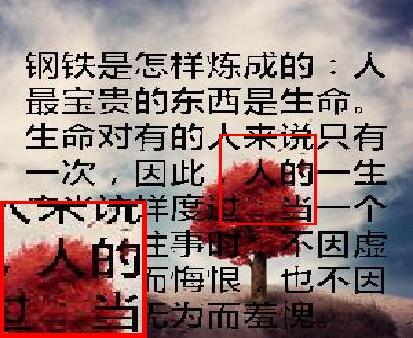}
}
\subfigure[MTNN:27.45dB.]{
\centering
\includegraphics[width=1.5in]{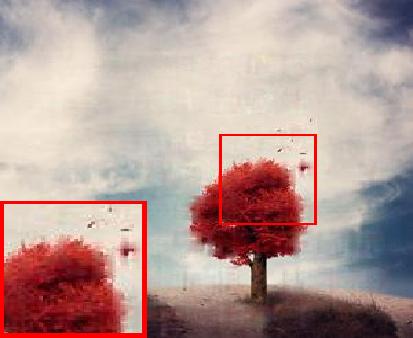}
}
\quad
\subfigure[T-TNN:27.81dB.]{
\centering
\includegraphics[width=1.5in]{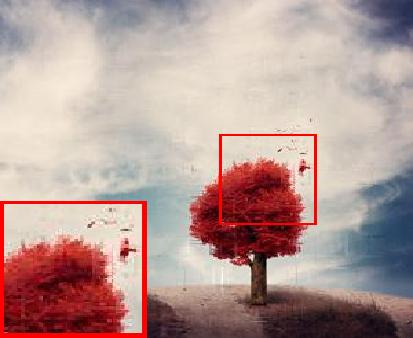}
}
\subfigure[T-TNNS:29.63dB.]{
\includegraphics[width=1.5in]{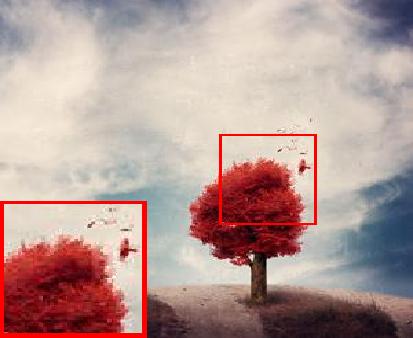}
}
\subfigure[SRTD:30.03dB.]{
\includegraphics[width=1.5in]{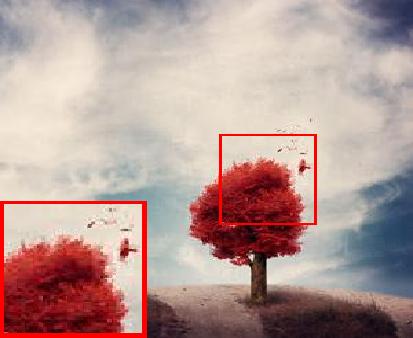}
}
\centering
\caption{The completion results with text mask.}\label{Fig:7}
\end{figure*}

\begin{figure*}[htbp]
\centering
\subfigure[Original image.]{
\centering
\includegraphics[width=1.5in]{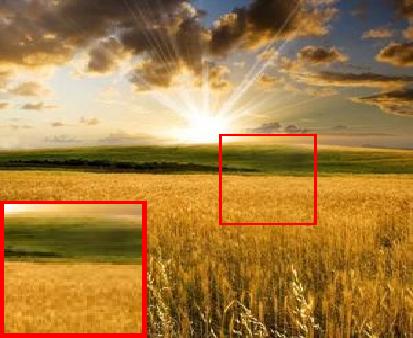}
}
\subfigure[text mask image.]{
\centering
\includegraphics[width=1.5in]{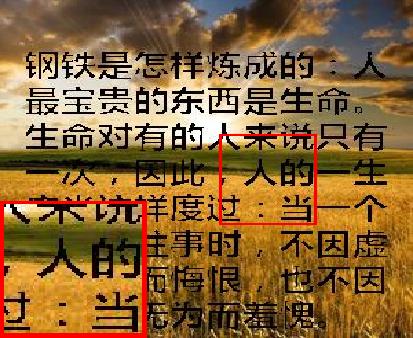}
}
\subfigure[MTNN:22.23dB.]{
\centering
\includegraphics[width=1.5in]{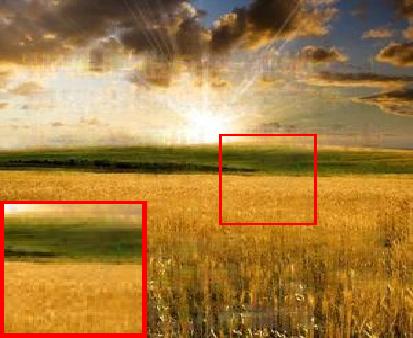}
}
\quad
\subfigure[T-TNN:23.42dB.]{
\centering
\includegraphics[width=1.5in]{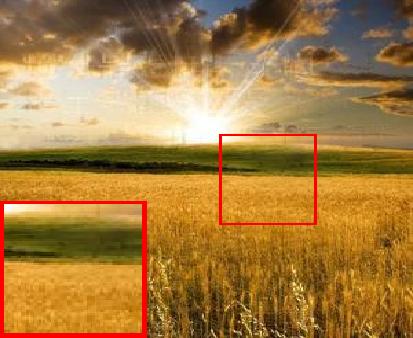}
}
\subfigure[T-TNNS:23.68dB.]{
\includegraphics[width=1.5in]{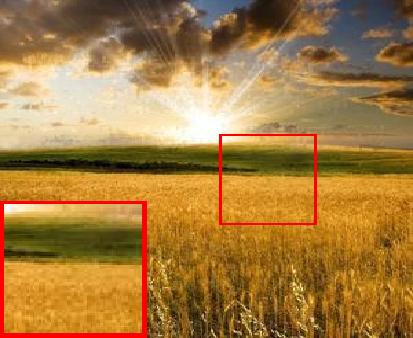}
}
\subfigure[SRTD:24.21dB.]{
\includegraphics[width=1.5in]{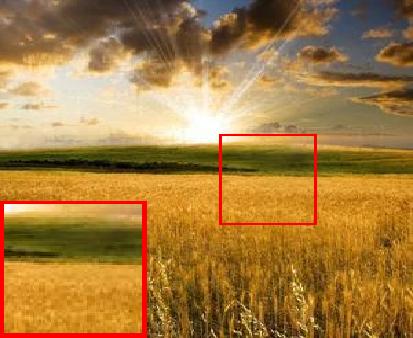}
}
\centering
\caption{The completion results with text mask.}\label{Fig:8}
\end{figure*}

In this part, we consider images which are corrupted by a text mask. It is a difficult task to remove the text, since the text is not randomly distributed in the image and it may cover some very important texture information. The text removal experiments results are shown in
Figure \ref{Fig:7} and Figure \ref{Fig:8}.

It can be seen from the recovered images that the proposed algorithm can recover the missing pixels  by the text mask noise very well. Moreover, we can see that the PSNR of the proposed SRTD is higher than that of MTNN, T-TNN and T-TNNS. Specifically, for the first image, the PSNR are 27.45, 27.81, 29.63 and 30.03 respectively for MTNN, T-TNN, and T-TNNS. For the second image, the PSNR are 22.23, 23.42, 23.68 and 24.21 respectively for MTNN, T-TNN, and T-TNNS. From both the PSNR values and visual effect, it can be demonstrated that the proposed method has a better performance.

\subsection{Video recovery with random mask}
\begin{figure*}[htbp]
\centering
\subfigure[Original 20th video.]{
\centering
\includegraphics[width=1.9in]{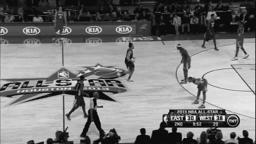}
}
\quad
\subfigure[35 $\%$ SR.]{
\centering
\includegraphics[width=1.9in]{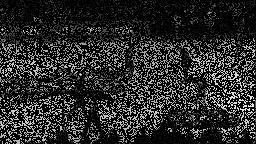}
}
\quad
\quad
\quad
\subfigure[T-TNN PSNR:23.12dB.]{
\centering
\includegraphics[width=1.9in]{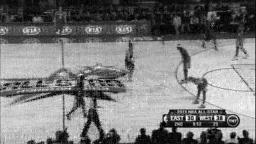}
}
\subfigure[T-TNNS PSNR:23.30dB.]{
\includegraphics[width=1.9in]{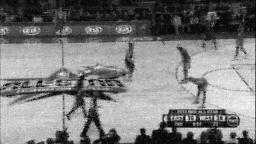}
}
\subfigure[SRTD PSNR:24.49dB.]{
\includegraphics[width=1.9in]{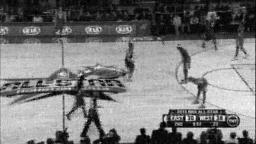}
}
\centering
\caption{The 20th frame of the basket video recovered by T-TNN, T-TNNS and SRTD at SR = 35\%.}\label{Fig:9}
\end{figure*}

\begin{figure*}[htbp]
\centering
\subfigure[Original 20th video.]{
\centering
\includegraphics[width=1.9in]{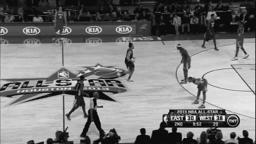}
}
\quad
\subfigure[25$\%$ SR.]{
\centering
\includegraphics[width=1.9in]{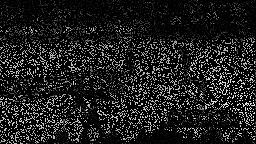}
}
\quad
\quad
\quad
\subfigure[T-TNN PSNR:21.65dB.]{
\centering
\includegraphics[width=1.9in]{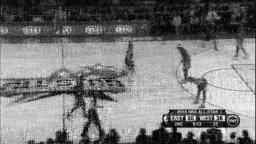}
}
\subfigure[T-TNNS PSNR:21.90dB.]{
\includegraphics[width=1.9in]{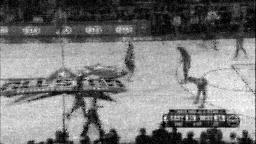}
}
\subfigure[SRTD PSNR:23.19dB.]{
\includegraphics[width=1.9in]{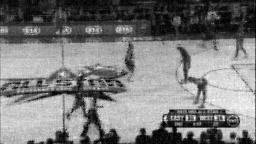}
}
\centering
\caption{The 20th frame of the basket video recovered by T-TNN, T-TNNS and SRTD at SR = 25\%.}\label{Fig:10}
\end{figure*}

Here we choose a gray basketball video which can been seen as a 3D tensor from YouTube.com with size $144\times256\times40$. The first two modes of the video correspond to the spatial variety, and the last mode corresponds to time changes. We set the SRs at 35$\%$ and 25$\%$ respectively. We compare the recovered PSNR value of proposed SRTD with T-TNN and T-TNNS, and show the contrast results of the 20th frame in Figure \ref{Fig:9} and Figure \ref{Fig:10}. Again, we observe that the performance of proposed SRTD is better than T-TNN and T-TNNS in PSNR values and visual effect.

\section{Conclusion}\label{sec5}
In this paper, we proposed a tensor completion approach SRTD based on the low-rank and sparse prior. In detail,
we used the tensor truncated nuclear norm based on T-SVD rather than the tensor nuclear norm used in most of the  existing methods, which can be regarded as a direct extension of matrix truncated nuclear norm. $\ell_{1}$-norm is used to describe the sparse prior of the tensor in a DCT domain, which is a general way to model the sparse property of tensors. A constrained optimization problem is formulated and then solved by ADMM iteration scheme. Experimental results showed that the proposed SRTD method performs better than MTNN, T-TNN and T-TNNS.

\section*{reference}
\bibliographystyle{elsarticle-num}

\end{document}